\documentstyle[amscd,amssymb,verbatim,12pt]{amsart}

\theoremstyle{plain}
\newtheorem{thm}{Theorem}
\newtheorem{lem}[thm]{Lemma}

\theoremstyle{definition}

\renewcommand{\rm}{\normalshape}

\newcommand{\epten}{\tilde{\otimes}_\epsilon}
\newcommand{\ep}{\epsilon}
\newcommand{\la}{\langle}
\newcommand{\ra}{\rangle}
\newcommand{\LR}{\Longleftrightarrow}
\newcommand{\N}{{\Bbb N}}
\newcommand{\R}{{\Bbb R}}

\begin{document}


\title{Property $\text{(u)}$ in $JH\epten JH$}
\author{Denny H.\ Leung}
\address{Department of Mathematics\\
         National University of Singapore\\
         Singapore 119260}
\email{matlhh@@leonis.nus.sg}



\maketitle

\begin{abstract}
It is shown that the tensor product $JH\epten JH$ fails Pe\l
cz\'{n}yski's  property (u).  The proof uses a result of Kwapie\'{n}
and Pe\l cz\'{n}yski on the main triangle projection in matrix spaces.
\end{abstract}

The Banach space $JH$ constructed by Hagler \cite{H} has a number of
interesting properties. For instance, it is known that $JH$ contains
no isomorph of $\ell^1$, and has property (S): every normalized weakly
null sequence has a subsequence equivalent to the $c_0$-basis. This
easily implies that $JH$ is $c_0$-saturated, i.e., every infinite
dimensional closed subspace contains an isomorph of $c_0$. In answer
to a question raised originally in \cite{H}, Knaust and Odell
\cite{KO} showed that every Banach space which has property (S) also
has Pe\l czy\'{n}ski's property (u).
In \cite{L}, the author showed that the Banach space $JH\epten JH$ is
$c_0$-saturated. It is thus natural to ask whether $JH\epten JH$ has
also the related properties (S) and/or (u).  In this note, we show
that $JH\epten JH$ fails property (u) (and hence property (S) as
well). Our proof makes use of a result, due to Kwapie\'{n} and 
Pe\l czy\'{n}ski, that the main triangle projection is unbounded in
certain matrix spaces.

We use standard Banach space notation as may be found in \cite{LT}.
Recall that a series $\sum x_n$ in a Banach space $E$ is called {\em
weakly 
unconditionally Cauchy}\/ (wuC) if there is a constant $K < \infty$
such that $\|\sum^k_{n=1}\epsilon_nx_n\| \leq K$
for all choices of signs $\epsilon_n = \pm 1$ and all $k \in \N$.  
A Banach space $E$ has {\em property}\/ (u) if
whenever $(x_n)$ is a weakly Cauchy sequence in $E$, there is a wuC
series $\sum y_k$ in $E$ such that $x_n - \sum^n_{k=1}y_k \to 0$
weakly as $n \to \infty$.
If $E$ and $F$ are Banach spaces, and $L(E',F)$ is the space of all
bounded linear operators from $E'$ into $F$ endowed with the operator
norm, then the tensor product $E\epten F$ is the closed subspace of
$L(E',F)$ generated by the weak*-weakly continuous operators of finite
rank. In particular, for any $x \in E$, and $y \in F$, one obtains an
element $x\otimes y \in E\epten F$ defined by $(x\otimes y)x' = \la x,
x'\ra y$ for all $x' \in E'$.

Let us also recall the definition of the space $JH$, as well as fix
some 
terms and notation.  Let $T = \cup^\infty_{n=0}\{0,1\}^n$ be the dyadic
tree. The elements of $T$ are called {\em nodes}.  If $\phi$ is a node
of the form $(\ep_i)^n_{i=1}$, we say that $\phi$ has {\em length}\/
$n$ and write $|\phi| = n$.  The length of the empty node is defined
to be $0$.  For $\phi, \psi \in T$ with $\phi = (\ep_i)^n_{i=1}$ and
$\psi = (\delta_i)^m_{i=1}$, we say that $\phi \leq \psi$ if $n \leq
m$ and $\ep_i = \delta_i$ for $1 \leq i \leq n$.  The empty node is
$\leq \phi$ for all $\phi \in T$.   Two nodes $\phi$ and $\psi$ are
{\em incomparable}\/ if neither $\phi \leq \psi$ nor $\psi \leq \phi$
hold. If $\phi \leq \psi$, we say that
 $\psi$ is a {\em
descendant}\/ of $\phi$, and we set
\[ S(\phi,\psi) = \{\xi: \phi \leq \xi \leq \psi\}. \]  
A set of the form $S(\phi,\psi)$ is called a {\em segment}, or more
specifically, an $m$-$n$ {\em segment}\/ provided $|\phi| = m,$ and
$|\psi| = 
n$. A {\em branch}\/ is a maximal totally ordered subset of $T$.
The set of all branches is denoted by $\Gamma$.  A branch $\gamma$
(respectively, a segment $S$) is
said to {\em pass through}\/ a node $\phi$ if $\phi \in \gamma$
(respectively, $\phi \in S$). If
$x: T \to \R$\ is
a finitely supported function and $S$ is a segment, we define (with
slight abuse of notation) $Sx = \sum_{\phi\in S}x(\phi)$.  In case $S
= \{\phi\}$ is a singleton, we write simply $\phi x$ for $Sx$. Similarly,
if $\gamma \in \Gamma$, we define $\gamma(x) = \sum_{\phi\in
\gamma}x(\phi)$.  A set of segments $\{S_1,\dotsc,S_r\}$ is {\em
admissible}\/ if they are pairwise disjoint, and there are $m, n \in
\N\cup\{0\}$ such that each $S_i$ is an $m$-$n$ segment.  The James
Hagler
space $JH$ is defined as the completion of the set of all finitely
supported functions $x: T \to \R$\ under the norm:
\[ \|x\| = \sup\left\{\sum^r_{i=1}|S_ix| : S_1,\dotsc ,S_r \mbox{ is an
admissible set of segments}\right\}. \]
Clearly, all $S$ and $\gamma$ extend to norm $1$ functionals on $JH$.
It is known that the set $T$ of all node functionals, and the set
$\Gamma$ of all branch functionals together span a dense subspace of
$JH'$ (cf.\ p.\ 301 of \cite{H}).
Finally, if $x: T \to \R$\ is finitely supported, and $n \geq 0$, let
$P_nx: T \to \R$\ be defined by 
\[ (P_nx)(\phi) = 
 \begin{cases}
   x(\phi)& \text{if $|\phi| \geq n$} \\
   0& \mbox{otherwise.}
 \end{cases} \]
Obviously, $P_n$ extends uniquely to a norm $1$ projection on $JH$,
which we denote again by $P_n$.  The proof of the following lemma is
left to the reader. We thank the referee for the succinct formulation.
 
\begin{lem} \label{proj}
For any $n \in \N$, construct a sequence $(\pi(1),\pi(2),\dotsc,\pi(n))$
by writing the odd integers in the set $\{1,\dotsc,n\}$ in increasing
order, followed by the even integers in decreasing order. Then 
\[
 (-1)^{\min(\pi(i),\pi(j))+1} = 1 \quad \LR \quad i + j \leq n+1.  
\]
\end{lem}

For any $n \in \N$ and $n \times n$ real matrix $M =
[M(i,j)]^n_{i,j=1}$, let $E(M)$ be the matrix
$[(-1)^{\min(i,j)+1}M(i,j)]$. Denote by $\sigma(M)$ the norm of $M$
considered as a linear map from $\ell^\infty(n)$ into $\ell^1(n)$,
i.e.,
\[ \sigma(M) = \sup\left\{ \sum^n_{i,j=1}a_ib_jM(i,j) : \sup_{1\leq
i,j\leq n}\{|a_i|, |b_j|\} \leq 1\right\}. \]

\begin{lem} \label{est}
There is a constant $C > 0$ such that for every $n \in \N$, there is
an $n \times n$ real matrix $M_n$ such that $\sigma(M_n) = 1$ and 
$\sigma(E(M_n)) \geq C \log n$.
\end{lem}

\begin{pf}
It follows easily from \cite[Proposition 1.2]{KP} that there are a
constant $C > 0$ and real $n \times n$ matrices $N_n = [N_n(i,j)]$ 
for every $n$
such that $\sigma(N_n) = 1$, and $\sigma([\epsilon(i,j)N_n(i,j)]) \geq C
\log n$, where 
\begin{equation*}
\epsilon(i,j) =
 \begin{cases}
   1& \text{if $i + j \leq n + 1$},\\
   -1& \text{otherwise}.
 \end{cases}
\end{equation*}
Let $\pi$ be the permutation in Lemma \ref{proj}.
Define $M_n(i,j) = N_n(\pi^{-1}(i), \pi^{-1}(j))$, $1 \leq i, j
\leq n$,
and let $M_n = [M_n(i,j)]$.  Clearly $\sigma(M_n) = \sigma(N_n)
= 1$ for all $n$.  Also,
\begin{multline*}
\sigma(E(M_n)) =
\sigma\left(\left[(-1)^{\min(\pi(i),\pi(j))+1}M_n(\pi(i),\pi(j))\right]\right)\\
 = \sigma\left(\left[\epsilon(i,j)N_n(i,j)\right]\right)
\geq C \log n,
\end{multline*}
as required.
\end{pf}

Let $\psi_1$ denote the node $(0)$, and 
$\psi_n = ({\displaystyle \overbrace{1\dots
1}^{n-1}0})$ for $n \geq 2$. For convenience, define $s_0 = 0$ and
$s_k = \sum^k_{i=1}i$ for $k \geq 1$.  Now choose a strictly
increasing sequence $(n_k)$ in $\N$, and a sequence of pairwise
distinct nodes
$(\phi_i)$ such that $\phi_i$ is a descendant of $\psi_k$ having
length $n_k$ whenever $s_{k-1} < i \leq s_k$, $k \in \N$.  For any $i
\in \N$, choose a branch $\gamma_i$ which passes through
$\phi_i$. If $i \leq s_k$, denote by $\phi(i,k)$ the node of length
$n_k$ which belongs to $\gamma_i$.  Finally, let $R_k =
[R_k(i,j)]^{s_k}_{i,j=s_{k-1}+1}$ be $k \times k$ real
matrices such that $\sum_k\sigma(R_k) < \infty$.
Then define a sequence of elements in $JH\epten JH$ as follows:
\[ U_l =
\sum^l_{k=1}\sum^{s_k}_{i,j=s_{k-1}+1}R_k(i,j)e_{\phi(i,l)}\otimes
e_{\phi(j,l)} 
\]
for $l \in \N$. Here, $e_\phi \in JH$ is the characteristic function
of the singleton set $\{\phi\}$.
Since the sequence
$(e_{\phi(i,l)})^{s_k}_{i=s_{k-1}+1}$ is isometrically equivalent to the
$\ell^1(k)$-basis whenever $k \leq l$,
\[ \left\|\sum^{s_k}_{i,j=s_{k-1}+1}R_k(i,j)e_{\phi(i,l)}\otimes
e_{\phi(j,l)}\right\| = \sigma(R_k), \]
and thus $\|U_l\| \leq
\sum_k\sigma(R_k) < \infty$ for any $l$. 

\begin{lem}
The sequence $(U_l)$ is a weakly Cauchy sequence in $JH\epten JH$.
\end{lem}

\begin{pf}
It is well known that a bounded 
sequence $(W_n)$ in a tensor product $E\epten
F$ is weakly Cauchy if and only if $(W_nx')$ is weakly Cauchy in $F$
for all $x' \in E'$.
Since $(U_l)$ is a bounded sequence, and $[T \cup \Gamma] = JH'$,
it suffices to show that $(U_lx')$
is weakly Cauchy in $JH$ for every $x'$ in  $T \cup
\Gamma$. Now for all $\phi \in T$, we clearly have $U_l\phi =
0$ for all large enough $l$. Next, consider any $\gamma \in
\Gamma$. If $\gamma$ does not pass through any $\psi_k$, then it
cannot pass through any $\phi(i,k)$ either. So $U_l\gamma = 0$ for all
$l$.  Otherwise, due to the pairwise incomparability of $(\psi_k)$,
there is a unique $k_0$ such that $\psi_{k_0} \in \gamma$.
If $\gamma$ is distinct from $\gamma_i$ for all $s_{k_0-1} < i \leq
s_{k_0}$, then again $U_l\gamma = 0$ for all sufficiently large
$l$. Now suppose $\gamma = \gamma_{i_0}$ , where $s_{k_0-1} < i_0 \leq
s_{k_0}$.  Then, for $l \geq k_0$,
\[ U_l\gamma = \sum^{s_{k_0}}_{j=s_{k_0-1}+1}R_{k_0}(i_0,j)e_{\phi(j,l)} .\]
Since each sequence $(e_{\phi(j,l)})^\infty_{l=k_0}$ is weakly Cauchy in
$JH$, so is $(U_l\gamma)$.
\end{pf}

Now if $JH\epten JH$ has property (u), then it is easy to observe that
there must be a block sequence of convex combinations $(V_r)$ of
$(U_l)$  such that $\sum (V_r - V_{r+1})$ is a wuC series. 
Write $V_r = \sum^{l_r}_{l=l_{r-1}+1}a_lU_l$ (convex combination),
where $(l_r)$ is a strictly increasing sequence in $\N$.
Fix $r \in \N$. For $s_{r-1}< i \leq s_r$, let $\xi_i$ be a branch
such that $\phi(i,l_{r+i-s_{r-1}})$ is the node of maximal length which
it shares with $\gamma_i$. Then if $s_{r-1} < i, j \leq s_r$, and $r
\leq l$,
\[ \la e_{\phi(i,l)}, \xi_j\ra = 1 \qquad \LR \qquad i = j \quad
\text{and} \quad l \leq l_{r+j-s_{r-1}}. \]
Hence, if $r \leq l$,
\[ \la U_l\xi_i, \xi_j\ra =
 \begin{cases}
  R_r(i,j)& \text{if $l \leq
       \min(l_{r+i-s_{r-1}},l_{r+j-s_{r-1}})$,}\\ 
  0& \text{otherwise.}
 \end{cases} \]
Thus, if  $s_{r-1} < i, j ,k \leq s_r$,
\[ \la V_{r+k-s_{r-1}}\xi_i, \xi_j\ra = 
 \begin{cases}
  R_r(i,j)& \text{if $k \leq \min(i,j)$,} \\
  0& \text{otherwise.}
 \end{cases} \]
It follows that 
\begin{multline*}
\left\la
\left\{\sum^{s_r}_{k=s_{r-1}+1}(-1)^{k+1-s_{r-1}}\left(V_{r+k-s_{r-1}}
- V_{r+k+1-s_{r-1}}\right)\right\}\xi_i, \xi_j\right\ra \\
= (-1)^{\min(i-s_{r-1},j-s_{r-1})+1}R_r(i,j). 
\end{multline*}
Notice that $k > s_{r-1}$ implies $l_{r+k-s_{r-1}} \geq l_{r+1} \geq
r$, hence  
\[ \la V_{r+k-s_{r-1}}P'_{n_r}\xi_i, P'_{n_r}\xi_j\ra = \la
   V_{r+k-s_{r-1}}\xi_i, \xi_j\ra. \]
Also, $(P'_{n_r}\xi_i)^{s_r}_{i=s_{r-1}+1}$ is isometrically
equivalent to the $\ell^\infty(r)$-basis.  Therefore,
\begin{multline*}
 \left\|\sum^{s_r}_{k=s_{r-1}+1}(-1)^{k+1-s_{r-1}}\left(V_{r+k-s_{r-1}}
- V_{r+k+1-s_{r-1}}\right)\right\| \\
\geq
\sigma\left(\left[(-1)^{\min(i-s_{r-1},j-s_{r-1})+1}R_r(i,j)\right]\right) =
\sigma(E(R_r)) . 
\end{multline*}
But since $\sum (V_r - V_{r+1})$ is a wuC series, there is a constant
$K < \infty$ (which may depend on the sequence $(R_k)$) such that 
\[ \left\|\sum^{s_r}_{k=s_{r-1}+1}(-1)^{k+1-s_{r-1}}\left(V_{r+k-s_{r-1}}
- V_{r+k+1-s_{r-1}}\right)\right\| \leq K \]
for any $r$. Consequently, $\sup_r\sigma(E(R_r)) \leq K$.

Now choose a strictly increasing sequence $(r_m)$ such that
$\lim_m 2^{-m}\log r_m = \infty$.  Then let
\[ R_k = 
 \begin{cases}
   \frac{M_{r_m}}{2^m}& \text{if $k = r_m$ for some $m$,}\\
   0& \text{otherwise,}
 \end{cases} \]
where $M_{r_m}$ is the matrix given by Lemma \ref{est}.
Then $\sum_k\sigma(R_k) = \sum_m2^{-m}\sigma(M_{r_m}) = 1$. So the
preceding argument yields a finite constant $K$ such that 
\[ K \geq \sup_m\frac{\sigma(E(M_{r_m}))}{2^m} \geq C 
\sup_m\frac{\log r_m}{2^m} ,
\]
contrary to the choice of $(r_m)$.  We have thus proved the following
result.

\begin{thm}
The Banach space $JH\epten JH$ fails property\/ {\rm (u)}.
\end{thm}



\begin{thebibliography}{99}

\bibitem{H} {\sc James Hagler},
 {\em A counterexample to several questions about Banach spaces},
 Studia Math.\ {\bf 60}(1977), 289-308.

\bibitem{KO} {\sc H.\ Knaust and E.\ Odell},
 {\em On $c_0$ sequences in Banach spaces},
 Israel J.\ Math.\ {\bf 67}(1989), 153-169.

\bibitem{KP} {\sc S.\ Kwapie\'{n} and A.\ Pe\l czy\'{n}ski},
 {\em The main triangle projection in matrix spaces and its 
  applications}, 
 Studia Math.\ {\bf 34}(1970), 43-68.

\bibitem{L} {\sc Denny H.\ Leung}, 
 {\em Some stability properties of $c_0$-saturated spaces},
 Math.\ Proc.\ Camb.\ Phil.\ Soc., {\bf 118}(1995), 287-301.

\bibitem{LT} {\sc Joram Lindenstrauss and Lior Tzafriri}, 
 Classical Banach Spaces I, Springer-Verlag, 1979.

\end{thebibliography}
\end{document}